\documentclass[11pt,a4paper]{article}
\usepackage{amssymb}

\usepackage[english]{babel}

\usepackage{amsmath}
\usepackage{amsthm}
\usepackage{verbatim}

\newtheorem{theorem}{Theorem}
\newtheorem{lemma}{Lemma}
\newtheorem{remark}{Remark}

\newtheorem{proposition}{Proposition}
\newtheorem{definition}{Definition}
\newtheorem{corollary}{Corollary\!\!}

\newcommand{\RR}{\mathbb{R}}
\newcommand{\CC}{\mathbb{C}}
\newcommand{\cH}{\mathcal{HM}}
\newcommand{\HCM}{\mathcal{HCM}}
\newcommand{\pcH}{\mathcal{{\widehat \cH}}}

\newcommand{\EE}{\mathbb{E}}
\newcommand{\cB}{\mathcal{B}}
\newcommand{\cS}{\mathcal{S}}

\newcommand{\tP}{{\hat P}}
\newcommand{\tH}{{\hat h}}

\newcommand{\cC}{\mathcal{C}}
\newcommand{\cCM}{\mathcal{CM}}

\newcommand{\cM}{\mathcal{M}}

\newcommand*\pFqskip{8mu}
\catcode`,\active
\newcommand*\pFq{\begingroup
        \catcode`\,\active
        \def ,{\mskip\pFqskip\relax}%
        \dopFq
}
\catcode`\,12
\def\dopFq#1#2#3#4#5{%
        {}_{#1}F_{#2}\biggl[\genfrac..{0pt}{}{#3}{#4};#5\biggr]%
        \endgroup
}

\def\a{{\alpha}}
\def\ii{{\rm i}}
\def\lacc{\left\{}
\def\lcr{\left[}
\def\lpa{\left(}

\def\racc{\right\}}
\def\rcr{\right]}
\def\rpa{\right)}

\def\Un{{\bf 1}}

\begin{document}

\title{Stieltjes functions of finite order and hyperbolic monotonicity}

\author{Lennart Bondesson\thanks{Department of Mathematics and Mathematical Statistics,  Ume\aa \, University, SE-90183 Ume\aa,  Sweden. $\qquad$ {\em Email:} \texttt{lennart.bondesson@umu.se}}$\qquad$ and $\quad$ Thomas Simon\thanks{Laboratoire Paul Painlev\'e, Universit\'e de Lille 1, Cit\'e Scientifique, F-59655 Villeneuve d'Ascq Cedex, France. $\qquad\qquad$
{\em Email:} \texttt{simon@math.univ-lille1.fr} }}


\maketitle

\parindent=0.5cm

\begin{abstract} A class of Stieltjes functions of finite type is introduced. These satisfy Widder's conditions on the successive derivatives up to some finite order, and are not necessarily smooth. We show that such functions have a unique integral representation, along some generic kernel which is a truncated Laurent series approximating the standard Stieltjes kernel. We then obtain a two-to-one correspondence, via the logarithmic derivative, between these functions and a subclass of hyperbolically monotone functions of finite type. This correspondence generalizes a representation of HCM functions in terms of two Stieltjes transforms earlier obtained by the first author.
\end{abstract}

2010 {\sl Mathematics Subject Classification.} 44A15, 60E10.\\

{\sl Keywords and phrases:}  Hyperbolic monotonicity, Stieltjes transform, Widder condition.

\section{Introduction and statement of the results}

\label{Sec:Intro}
This paper is devoted to certain subclasses of real functions defined on $(0,\infty).$ Unless otherwise explicitly stated, all functions will be supposed of this kind in the sequel.  A non-negative function $f$ is called a Stieltjes function ($f\in\cS$ for short) if there exists $a \ge 0$ and a non-negative measure $\mu(dt)$ on $[0,\infty)$ integrating $(1+ t)^{-1}$ such that
\begin{equation}
\label{SJ}
f(x)\; =\; a\; +\; \int_0^\infty \frac{1}{x+t}\;\mu(dt).
\end{equation}
Introduced by Stieltjes for the purposes of the moment problem, such functions are important for questions related to potential theory and infinite divisibility - see \cite{SSV} and the references therein for a recent account, among other topics. Notice that in many instances, it is useful to extend Stieltjes functions to the whole cut plane $\CC\backslash \RR^-$. In this paper however, we will stay within the realm of functions of one real variable.

It is plain by dominated convergence that a function $f\in\cS$ is also completely monotone ($f\in\cCM$ for short), in other words $f$ is smooth and
\begin{eqnarray}
\label{CM}
(-1)^{n}f^{(n)}\, \ge\, 0
\end{eqnarray}
for all $n\ge 0$ where, here and throughout, $f^{(n)}$ stands for the $n-$th derivative of $f.$ Recall from Bernstein's theorem - see e.g. Theorem 1.4 in \cite{SSV} - that $f\in\cCM$ if and only if there exists a non-negative measure $\mu (dt)$ on $[0,\infty)$ (the so-called Bernstein measure) such that
$$f(x)\; =\; \int_0^\infty e^{-xt}\, \mu(dt).$$
It is also easy to see that if $f\in\cS,$ then $x f$ is a Bernstein function ($f\in\cB$ for short), that is a non-negative function whose derivative is completely monotone - see again \cite{SSV} for an account. More generally, it was shown by Widder - see Theorem 10.1 in \cite{W0} - that $f\in\cS$ if and only if $f$ is smooth and
\begin{eqnarray}
\label{Sti0}
(x^n f)^{(n)}\,\in\,\cCM \quad \forall\,n\ge 0.
\end{eqnarray}
The proof of this result was recently simplified in \cite{S}, in the broader framework of generalized Stieltjes transforms. Another theorem by Widder - see Theorem 12.5 in \cite{W1} or Theorem 18b p. 366 in \cite{W} - states that a non-negative function $f$ is in $\cS$ if and only if it is smooth and such that
\begin{eqnarray}
\label{Sti1} (-1)^{n-1}(x^n f)^{(2n-1)}\; \ge\; 0\quad\forall\, n\ge 1.
\end{eqnarray}
Notice that the formal equivalence between (\ref{Sti0}) and (\ref{Sti1}), which is partly explained in Lemma 12.52 of \cite{W1}, is not immediate.

Finite type versions of (\ref{CM}) and (\ref{Sti0}) have been studied in the literature. Following \cite{Wi}, we will say that a function $f$ is $k-$monotone ($f\in\cM_k$ for short) for some $k\ge 2$ if it is in $\cC^{k-2}$ and such that $(-1)^{n}f^{(n)}$ is non-negative, non-increasing and convex for $n = 0,\ldots, k-2$. As in \cite{Wi}, we will say that $f\in\cM_1$ if $f$ is non-negative and non-increasing. These functions have been characterized in Theorem 1 in \cite{Wi} - see also Lemma 17.4.1 p. 306 in \cite{LB97}, which states that $f\in\cM_k$ if and only if there exists a non-negative measure $\mu_k (dt)$ such that
$$f(x)\; =\; \int_0^\infty \lpa 1-\frac{xt}{k}\rpa_{\! +}^{\! k-1}\!\! \mu_k(dt).$$
Notice that this result recovers Bernstein's theorem by letting $k\to \infty,$ identifying the exponential kernel $e^{-xt}$ at the limit of the integrand, and applying Helly's selection principle to the sequence $\{\mu_k\}$ - see the remark p. 310 in \cite{LB97} for details. More recently, in \cite{P}, it was shown that a non-negative function $f$ satisfies (\ref{Sti0}) for $n=1,\ldots, k$ if and only if it is in $\cCM$ and its Bernstein measure has a $k-$monotone density satisfying a certain integrability property - see Theorem 1.3 therein. This result also retrieves Theorem 10.1 in \cite{W0} by the same selection argument - see Corollary 1.5 therein.\\

In this paper, we will obtain a finite type version of (\ref{Sti1}). This motivates the following definition, which is inspired by \cite{W1} and \cite{Wi}. Here and throughout, we will consider derivatives in the measure sense. With this convention, the first derivative of a convex function on $(0,\infty)$ is its right derivative, whereas its second derivative is a non-negative measure.

\begin{definition}
\label{Stn}
For $k\ge 1,$ a non-negative function $f$ satisfying {\em (\ref{Sti1})} for $n = 1,\ldots, k$ is said to be a $k-$Stieltjes function ($f\in\cS_k$ for short).
\end{definition}

Notice that for $k\ge 2,$ a function $f\in\cS_k$ is of class $\cC^{2k-3}$ with $(-1)^{k-1} (x^k f)^{(2k-3)}$ a convex function. Observe also that $k$-Stieltjes functions need not even for $n=0$ satisfy (\ref{Sti0}), since the condition is on a finite number of derivatives only. On the other hand, it will be shown in Proposition \ref{SkMk} below that if $f\in\cS_k$ for some $k\ge 2,$ then $(xf)'\in\cM_{k-1}.$ It is clear that $\{\cS_k\}$ is a decreasing family with $\cS_k\downarrow\cS.$ To state our results, we need some further notation. Introduce the following family of alternating sign polynomials
$$P_k(x)\; =\; \sum_{n=0}^k \lpa \!\!\!\begin{array}{c} 2k\\ n+k \end{array}\!\!\!\rpa (-x)^{n}.$$
Observe that $P_k(x)$ is the polynomial part of the finite Laurent series $(1-x)^k(1-x^{-1})^k.$ Set $\tP_k = P_k -P_k(0)$ and consider the following family of non-negative kernels on $(0,\infty)\times[0,\infty):$
$$\Phi_k(x, t)\; =\; \frac{1}{x}\lpa  P_k(t x^{-1}) \Un_{\{x\ge t\}}\,-\, \tP_k(x t^{-1}) \Un_{\{x<t\}}\rpa.$$
Observe that $x\Phi_k(x, 0) = P_k(0) = \binom{2k}{k}$ and that $x\Phi_0(x, t) = \Un_{\{x\ge t\}}.$ The fact that the kernels $\Phi_k(x,t)$ are everywhere non-negative is a direct consequence of the decreasing character of the coefficients of $P_k,$ which implies $P_k(y)\ge 0$ and $\tP_k(y) \le 0$ for all $y\in[0,1].$ More generally, it will be proved in Proposition 1 below that for every $k\ge 2$ and $t > 0,$ the functions $x\mapsto \Phi_{k-1}(x,t)$ belong to $\cS_k.$ \\

Our first main result is the following characterization.

\begin{theorem}
\label{Widderlike}
For $k\ge 2,$ one has $f\in\cS_k$ if and only if there exists $a_k\ge 0$ and a non-negative measure $\mu_k(dt)$ on $[0,\infty)$ integrating $(1 + t)^{-1}$ such that
\begin{eqnarray}
\label{Wid}
f(x)\; =\; a_k\, +\, \int_0^\infty\! \Phi_{k-1} (x,t)\, \mu_k(dt).
\end{eqnarray}
\end{theorem}

Recall that in the case $k=1,$ the condition $f\in\cS_1$ means that $x f$ is non-decreasing and we hence have the obvious representation
$$f(x)\; =\; \frac{1}{x} \int_0^x \mu_1(dt)\; =\; \int_0^\infty\! \Phi_0 (x,t)\, \mu_1(dt)$$
for some non-negative measure $\mu_1$ on $[0,\infty)$ which, however, might not integrate $(1 + t)^{-1}.$ Observe also that
\begin{equation}
\label{Stic}
\binom {2k}{k}^{-1} \Phi_k(x, t)\; \longrightarrow\; \frac{1}{x+t}
\end{equation}
as $k\to \infty$ for all $x>0$ and $t\ge 0,$ so that again, applying Helly's selection principle one retrieves Widder's characterization of $\cS$ given in (\ref{Sti1}). It is plain that the sets $\cS_k$ are convex cones of functions, and the above result together with the argument of Proposition 1 in \cite{H} shows that they are closed with respect to pointwise limits. It is not clear whether these closed convex cones have abstract extensions leading to interesting invariant properties, as is the case for $\cS$ - see \cite{H, Be}. Recall that the extension of $\cS$ to $\CC\backslash \RR^-$ implies a complex inversion formula which is valid when $\mu$ in (\ref{SJ}) has a density, and which is well-known as the Perron-Stieltjes inversion formula - see e.g. Theorem 7b p. 340 in \cite{W}. In our finite type framework there is also an inversion formula for $\mu_k$, which has a real-variable character and is more directly connected to the kernels $\Phi_k$ and the conditions (\ref{Sti1}) - see Remark \ref{Inv} below. \\

The proof of the characterization of $\cS$ by the set of conditions (\ref{Sti1}), which is more or less the topic of the whole Chapter 8 in \cite{W}, is lenghty. It hinges on the construction of a certain jumping operator connected to the Perron-Stieltjes inversion formula. Our proof goes partly along Widder's lines, but the main difference is that it relies on the {\em non-smoothness} of the kernels $\Phi_k.$ More precisely, applying the $k-$th condition (\ref{Sti1}) to $\Phi_{k-1}$ yields a Dirac mass  - see Remark \ref{Expl}. This observation, which explains why the truncated Laurent series $\Phi_k$ are the relevant approximations of the Stieltjes kernel in our context, allows us to find the measure $\mu_k$ in a constructive way, starting from the convexity assumption in Definition 1 (b) and then integrating. The integration procedure works and gives the right growth order for $\mu_k,$ because the assumption $f\in\cS_k$ forces the function $f$ to have a certain boundary behaviour at zero - see Proposition \ref{SkMk}. Overall, the problem that we consider in this paper is more complicated than the problems in  \cite{P, S, W} because of its non-smooth character, and our arguments are also more intricate. \\

As mentioned before, Stieltjes functions appear in questions related to infinite divisibility. This is mainly due to the aforementioned property that $f\in\cS\Rightarrow xf\in\cB$ and we refer e.g. to Chapters 7 and 8 in \cite{SSV} for more on this topic. A further instance is the following notion, introduced by Thorin and the first author in the late 1970's: a function $f$ is said to be hyperbolically completely monotone ($f\in\HCM$ for short) if the function $f(uv)f(uv^{-1})$ is completely monotone in the variable $w = v +v^{-1},$ for every $u > 0.$ This apparently technical definition is actually quite robust, and a remarkable feature of the class $\HCM$ is that such functions appear both as Laplace transforms and densities of infinitely divisible distributions. It turns out that functions in $\HCM$ are pointwise limits of functions of the type
$$C x^{\beta-1}\prod_{i=1}^N (1+c_ix)^{-\gamma_i}$$
with all parameters positive except $\beta\in\RR.$  The connection with Stieltjes functions is obtained by the following representation, which is given as (5.2.3) in \cite{LB92} and is a consequence of Theorem 5.3.1 therein: one has $f\in\HCM$ if and only if
\begin{eqnarray}
\label{CanHCM}
f(x) \; =\; C x^{\beta-1}h_1(x)h_2(x^{-1}),
\end{eqnarray}
with $C \ge 0, \beta \in \RR,$ and $-(\log h_i)' \in \cS$ for $i=1,2.$ From this representation, it is clear that $f^p\in\HCM$ for all $p > 0$ if $f\in\HCM$ and that $fg\in\HCM$ whenever $f, g \in\HCM.$ A deeper property is that $\HCM$ is also stable by multiplicative convolution. We refer to Chapters 3-5 in \cite{LB92} for more details on this notion. \\

In this paper, we will obtain a finite type version of (\ref{CanHCM}). This motivates the following definition, which is rephrased from the main definition of \cite{LB97}.

\begin{definition}
\label{defHM}
A non-negative function $f$ is called $k-$hyperbolically monotone ($f\in\cH_k$ for short) if, for every $u>0,$ the function $f(uv)f(uv^{-1})$ is $k-$monotone in the variable $w = v +v^{-1}.$
\end{definition}

Again, we see that $\{\cH_k\}$ is a decreasing family with $\cH_k\downarrow\HCM.$
The Leibniz formula shows that $\cH_k$ is closed with respect to multiplication, and it is easy to see that it is also closed with respect to pointwise limits and to the transformation $f \to \tilde f(x)= f(x^{-1}).$ In \cite{LB97}, it was shown among other results that $\cH_k$ is closed with respect to multiplicative convolution. We also refer to \cite{BB} for further connections between the class $\cH_k$ and infinite divisibility.

In the case $k=1,$ it is not difficult to see that $f\in\cH_1$ if and only if
$$f(u_1v_1^{-1})f(u_2v_2^{-1}) \; \ge \; f(u_1v_2^{-1})f(u_2v_1^{-1})$$
for every $u_1 < u_2$ and $v_1 < v_2.$ This means that the kernel $f(xy^{-1})$ is ${\rm TP}_2$ on $(0,\infty)\times (0,\infty)$ or equivalently - see Theorem 4.1.8 in \cite{K},  that $y \mapsto f(e^y)$
is log-concave on its support which is necessarily a closed interval. Hence, there is a  canonical representation: one has $ f\in \cH_1$ if and only if
\begin{eqnarray}
\label{H1}
f(x)\;= \;C \exp\lcr -\int_{x_0}^x \frac{\psi(y)}{y}dy\rcr,
\end{eqnarray}
with $ C, x_0 > 0$ suitably chosen, and $\psi$ a non-decreasing function (possibly taking the values $\pm \infty$). Separating the positive and negative parts of $\psi,$ it is an easy exercise to transform this representation into
\begin{eqnarray}
\label{CanHM1}
f(x) \; =\; C x^{\beta-1}h_1(x)h_2(x^{-1})
\end{eqnarray}
with $C \ge 0, \beta \in \RR,$ and $-(\log h_i)' \in \cS_1$ (possibly taking the value $+\infty$) for $i=1,2.$ In particular, we see that $f^p \in\cH_1$ for every $p >0.$ In the case $k\ge 2$ however, the connection between $\cH_k$ and totally positive kernels of higher order is lost in general. Moreover, it is possible to exhibit functions $f\in\cH_k$ such that $f^p\not\in\cH_k$ for some $p >0$ - see Remarks 4(b) and 7(b). Having in mind an exponential representation of the type (\ref{CanHCM}) or (\ref{CanHM1}), it is hence natural to introduce the following definition:

\begin{definition}
\label{PR}
A non-negative function $f$ such that $f^p\in\cH_k$ for every $p>0$
is called power regular $\cH_k$ ($f\in\pcH_k$ for short).
\end{definition}

Our second main result is the following characterization.

\begin{theorem}
\label{Bondlike}
For every $k\ge 2,$ one has $f\in\pcH_k$ if and only if
$$f(x)\; =\; C x^{\beta-1}h_1(x)h_2(x^{-1})$$
with $C \ge 0, \beta \in \RR,$ and $-(\log h_i)' \in \cS_k$ for $i=1,2.$
\end{theorem}

This result gives a constructive procedure to find functions in the set $\pcH_k$, starting either from two functions satisfying (\ref{Sti1}) for every $n\le k$ or, by Theorem 1, from two non-negative reals and two non-negative measures on $(0,\infty)$ integrating $(1+t)^{-1}.$ An example of the latter construction is provided in Section 4.2. It remains an open problem to find a canonical representation for all $k-$hyperbolically monotone functions. This problem seems however uneasy because $\cH_k$ is not closed with respect to positive powers, which shows that the canonical representation, if any, should not have an exponential type.

The closedness of $\pcH_k$ with respect to positive powers plays a crucial role in our argument because it allows to linearize the problem - see Lemma \ref{HMD}, making any function in $\pcH_{k+1}$ in one-to-one correspondence with a parametrized set of $k-$monotone functions in the hyperbolic variable $v+v^{-1}.$ The remainder of the proof is then an analysis on this set of functions, whose initial conditions establish the connection with $\cS_k$ by a Taylor expansion. An
unexpected feature, which is a consequence of both Theorem 1 and the specific nature of our kernels $\Phi_k(x,t)$, is that these initial conditions determine the whole $k-$monotonicity property of these functions - see Remark \ref{Ini}. Let us finish this introduction with the following further characterization of $\pcH_k$, which is a simple consequence of Theorem 1, Theorem 2, and Lemma \ref{Form} below, and which we state without proof. An example of this characterization is provided in Section 4.1.

\begin{corollary}
\label{Bondesson}
For every $k\ge 2,$ one has $f\in\pcH_k$ if and only if
$$(-1)^j (x^j (\log f)')^{(2j-1)}\; \ge\; 0$$
for every $j\le k.$
\end{corollary}

The paper is organized as follows. In Section 2 we prove Theorem 1, and in Section 3 we prove Theorem 2. In Section 4 we consider some explicit interesting examples, whereas the Appendix is devoted to a technical and rather surprising Lemma related to the if part of Theorem 2.

\section{Proof of Theorem \ref{Widderlike}}

\subsection{Proof of the if part}

We will first investigate regularity properties of the kernel $\Phi_k(x,t),$ which are less immediate than those of the finite type kernels involved in \cite{Wi} and \cite{P}. For symmetry reasons, it will be more pleasant to consider the kernel
\begin{equation}
\label{Msika}
\Psi_k(x,t)\; = \; x \Phi_k(x,t)\; = \psi_k (x t^{-1}),
\end{equation}
where we have set $\psi_k(y) = \Psi_k (y, 1) =  P_k(y^{-1})\Un_{\{y \ge 1\}}-\tP_k(y)\Un_{\{y < 1\}}.$ The kernel $\Psi_k$ and the function $\psi_k$ will play some role in other parts of the paper. Our analysis relies on two lemmas which will be useful in other parts of the paper as well. The first one is obtained by an elementary induction, starting from the last derivative.

\begin{lemma}
\label{Poly}
Let $P$ be a real polynomial of degree $l.$ Suppose $(-1)^i P^{(i)}(0) >  0$ and $(-1)^i P^{(i)}(1) >  0$ for all $i=0,\ldots, l.$ Then $P$ is $l-$monotone on $[0,1].$ In particular, $P$ is positive on $[0,1].$
\end{lemma}
The second one has a symmetry character and is reminiscent of Lemma 3.11 in \cite{W1}. It consists in two identities between differential operators which are easily checked on polynomials, and thus on all functions by an identification of the coefficients. Alternatively, these identities can be obtained from the Leibniz formula.  Introduce the linear differential operator
$$\Theta_n(h) \; =\; x^n (x^{n-1} h)^{(2n-1)}$$
acting on any function $h$ which is regular enough. Let further $\tH(x) = -h(x^{-1}).$

\begin{lemma}
\label{Form}
For any $h$ regular enough, one has
$$ \Theta_n(h)(x)\; =\; (x^{2n-1} h^{(n)}(x))^{(n-1)}\; =\;\Theta_n(\tH)(x^{-1}).$$
\end{lemma}

\medskip

We can now state the main result of this paragraph.

\begin{proposition}
\label{IP1} For every $k\ge 1, t > 0,$ the function $x\mapsto \Phi_k(x,t)$ belongs to $\cS_{k+1}.$
\end{proposition}

\proof  We begin with an analysis of the function $\psi_k,$ which is smooth on $(0, \infty)$ except possibly at $y =1.$ Evaluating
$$P_k (1)\; =\; \lpa \!\!\!\begin{array}{c} 2k\\ k \end{array}\!\!\!\rpa\;\sum_{n=0}^k \frac{(-k)_n}{(k+1)_n}\; =\; \lpa \!\!\!\begin{array}{c} 2k\\ k \end{array}\!\!\!\rpa\;\pFq{2}{1}{-k,1}{k+1}{1}\; =\; \frac{1}{2} \lpa \!\!\!\begin{array}{c} 2k\\ k \end{array}\!\!\!\rpa\; =\;-\tP_k(1),$$
where we have used the standard notation for Pochhammer symbols and the hypergeometric function, and applied the Chu-Vandermonde identity  - see e.g. Corollary 2.2.3 in \cite{AAR}, shows that $\psi_k$ is continuous at $y =1.$ Similarly, we compute
$$P_k^{(i)} (1)\; =\;\tP_k^{(i)} (1)\; =\; (-1)^i \, i!\lpa \!\!\!\begin{array}{c} 2k\\ k +i \end{array}\!\!\!\rpa\;\pFq{2}{1}{i-k,i+1}{k+i+1}{1}\; =\; (-1)^i\,i!\lpa \!\!\!\begin{array}{c} 2k-i-1\\ k-1 \end{array}\!\!\!\rpa$$
for all $i = 1,\ldots, k.$ On the other hand, setting $Q_k (y) = \tP_k(y^{-1}),$ we obtain after some analogous computations
$$Q_k^{(i)} (1)\; =\;(-1)^{i+1} \, i!\lpa \!\!\!\begin{array}{c} 2k\\ k +1 \end{array}\!\!\!\rpa\;\pFq{2}{1}{1-k,i+1}{k+2}{1}\; =\;(-1)^{i+1} \, i!\lpa\frac{(k-i+1)_{k-1}}{(k-1)!}\rpa$$
for all $i \ge 1.$ For $i=1,\ldots, k,$ we get
$$Q_k^{(i)} (1) \;=\; (-1)^{i+1}\,i!\lpa \!\!\!\begin{array}{c} 2k-i-1\\ k-1 \end{array}\!\!\!\rpa\; =\;- \tP_k^{(i)} (1),$$
whereas for $i =k+1,\ldots,2k-1,$ we have $Q_k^{(i)} (1) = 0 = \tP_k^{(i)} (1).$ All of this shows that $\psi_k$ is of class $\cC^{2k-1}$ at $y =1$, and the function $x\mapsto \Phi_{k} (x,t)$ is hence $\cC^{2k-1}$ on $(0,\infty),$ too.\\

We next prove that
\begin{eqnarray}
\label{Sti2} \Psi_{n,k} (x)\; =\;  (-1)^{n-1} (x^{n-1} \psi_k(x))^{(2n-1)}\; \ge\; 0
\end{eqnarray}
for all $n =1,\ldots, k$ and $x >0.$ Suppose first $x < 1.$ Then, since
$$x^{n-1}\psi_k(x)\; =\; - x^{n-1} \,\tP_k(x)\; =\; - \binom{2k}{k} \sum_{i=0}^{k-1} \frac{(-k)_{i+1}}{(k+1)_{i+1}}\, x^{i+n},$$
we obtain
$$\Psi_{n,k}(x)\; =\; (2n-1)!\,\sum_{i=0}^{k-n} \lpa \!\!\!\begin{array}{c} 2k\\ k-n-i \end{array}\!\!\!\rpa \lpa \!\!\!\begin{array}{c} 2n+i-1\\ 2n-1 \end{array}\!\!\!\rpa (-x)^i.$$
For every $i =0,\ldots, k-n$ we next evaluate
$$\frac{(-1)^i \Psi_{n,k}^{(i)}(0)}{(2n-1)!i!}\, =\,   \lpa \!\!\!\begin{array}{c} 2k\\ k-n-i \end{array}\!\!\!\rpa \lpa \!\!\!\begin{array}{c} 2n+i-1\\ 2n-1 \end{array}\!\!\!\rpa \;\;\mbox{and}\;\; \frac{(-1)^i \Psi_{n,k}^{(i)}(1)}{(2n-1)!i!}\, =\,   \lpa \!\!\!\begin{array}{c} 2k-2n-i\\ k-n-i \end{array}\!\!\!\rpa \lpa \!\!\!\begin{array}{c} 2n+i-1\\ 2n-1 \end{array}\!\!\!\rpa,$$
where the first identity is immediate and the second one is obtained by the Chu-Vandermonde identity. Lemma \ref{Poly} implies then that (\ref{Sti2}) holds for all  $n=1,\ldots, k$ and $x < 1.$ When $x >1,$ we compute
$$\Psi_{n,k} (x) =  (-1)^{n-1} (x^{n-1} (\psi_k(x) -P_k(0)))^{(2n-1)}=  (-1)^{n-1} (x^{n-1} Q_k(x))^{(2n-1)} = x^{-2n} \Psi_{n,k} (x^{-1}),$$
where the first equality is obvious and the third one follows from Lemma \ref{Form}.  Putting everything together implies that (\ref{Sti2}) holds for every $n =1,\ldots, k$ and all $x \in(0,\infty).$ In particular, the function $x\mapsto \Phi_k(x,t)$ satisfies (\ref{Sti1}) for every $n =1,\ldots, k.$ \\

We finally consider the continuous function
$$\pi_k(x) \; =\; (-1)^{k} (x^k \psi_k(x))^{(2k-1)}\; =\; -x\Psi_{k,k}(x)\; +\; (-1)^k (2k-1) (x^{k-1} \psi_k(x))^{(2k-2)},$$
where in the second equality we have used the chain rule. If $x < 1,$ we have $\Psi_{k,k}(x) = (2k-1)!$ and $\pi_k$ has hence constant derivative $-(2k)!.$ If $x >1,$ Lemma \ref{Form} (or a direct computation) shows that $\Psi_{k,k}(x) = (2k-1)!x^{-2k}$ and this easily implies that $\pi_k$ has zero derivative. In particular, the function $\pi_k$ is convex, which means that the function $x\mapsto \Phi_k(x,t)$ fulfils the required convexity property for $\cS_{k+1}.$

\qed

\bigskip

It is easy to see that Proposition \ref{IP1} proves the if part of Theorem 1. Suppose indeed that $f$ has the representation (\ref{Wid}) for some $k \ge 2.$ Since $a + g\in\cS_k$ whenever $g\in\cS_k$ for every $a\ge 0,$ it suffices to consider the case $a_k = 0.$ Moreover, at each $x > 0$, it is plain by definition that the $i-$th derivative of $x\mapsto \Phi_{k-1}(x,t)$ is bounded in $t$ by $K_{i,x} (t^{i}\wedge t^{-i})$ for some finite constant $K_{i,x}.$ Since $\mu_k$ integrates all such functions by assumption, we can apply the dominated convergence theorem, and the linearity of the integral for the convexity property, to conclude by Proposition \ref{IP1} that $f$ belongs to $\cS_k.$

\qed

\medskip

\begin{remark}
\label{Expl}
{\em The end of the proof of Proposition \ref{IP1} shows that
$$\pi_k^{(2)}\; =\; (2k)!\, \delta_1$$
where, here and throughout, we have set the standard notation $\delta_a$ for a Dirac mass at $a$. This property allows to construct the kernel $\Phi_k$ in a recursive way, by successive integration choosing the appropriate boundary terms in order to ensure the required regularity and so that (\ref{Sti1}) holds at all intermediate levels. In this respect, the functions $x\mapsto\Phi_k(x,t)$ can be viewed as the ``fundamental solutions" of (\ref{Sti1}) up to order $k+1$ since we have
$$(-1)^k (x^{k+1} \Phi_k(x,t))^{(2k+1)}\; =\; (2k)!\, t^{-k}\, \delta_t.$$
We refer to Section 4.3 below for another recursive formula connecting $\Phi_k$ and $\Phi_{k+1}.$}
\end{remark}

\subsection{Proof of the only if part}

We begin with a proposition having an independent interest, and crucial for our purposes. For simplicity, we will set
$$g\; =\; x f\qquad\mbox{and}\qquad\varphi_n\; = \;x^{2n-1} g^{(n)}$$
for every suitable $n\ge 1.$

\begin{proposition}
\label{SkMk}
Suppose $f\in\cS_k. $ Then $g'\in \cM_{k-1}.$ Moreover, one has $(-1)^{k-1}\varphi_k^{(k-2)}(0+) \ge 0$ and $\varphi_k^{(j)}(0+) =0$ for all $j=0,\ldots, k-3.$
\end{proposition}

\proof  Observe first that by the first equality in Lemma \ref{Form}, if $f\in\cS_k$ then
$$(-1)^{n-1}\varphi_n^{(n-1)}\; \ge\;  0$$
for every $n =1,\ldots, k-1,$ and $(-1)^{k-1}\varphi_k^{(k-1)}$ is a non-negative measure on $(0,\infty).$ We now proceed by induction on $k.$

\medskip

If $k=2,$ then $\varphi_2'$ is a negative measure, which implies $\varphi_2(0+) > -\infty.$ Supposing $\varphi_2(0+) > 0,$ then recalling $\varphi_2 = x^3 g^{(2)},$ we see by integration that $g'(0+) =-\infty.$ This is a contradiction since $g$ is non-decreasing, by the fact that $f$ also belongs to $\cS_1$. Hence $\varphi_2(0+) \le 0$ and $g^{(2)}\le 0,$ which implies that $g'\in\cM_1$ as required.

\medskip

If $k=3,$ then $\varphi_3^{(2)}$ is a non-negative measure, which implies $\varphi_3'(0+) < +\infty$ and $\varphi_3(0+) > -\infty.$ If $\varphi_3(0+) < 0,$ then by integration one has $g^{(2)}(0+)=+\infty,$ a contradiction with the case $k=2.$ Moreover, the case $k=2$ also implies $\varphi_3  + 3x\varphi_2 \le  0$ and since $\varphi_2(0+) > -\infty,$ we must have $\varphi_3(0+) \le 0.$ This shows $\varphi_3(0+) = 0.$ Supposing now $\varphi_3'(0+) < 0,$ then again this implies the contradiction $g^{(2)}(0+)=+\infty.$ Hence we have $\varphi_3(0+) = 0, \varphi_3'(0+) \ge 0$ and $g^{(3)}\ge 0,$ which together with the case $k=2$ implies that $g'\in\cM_2.$

\medskip

We now set $k\ge 4$ and suppose that the property has been shown up to rank $k-1.$ Since $(-1)^{k-1} \varphi_k^{(k-1)}$ is a non-negative measure, a direct induction shows that all right derivatives $\varphi_k^{(j)}(0+)$ exist for $j=0,\ldots, k-2.$ A further induction based on L'H\^ospital's rule and the fact, given by the case $k=3,$ that $\varphi_3(0+) = 0,$ shows that $\varphi_k^{(j)}(0+) =0$ for $j=0,\ldots, k-3.$

Supposing $(-1)^{k-1} \varphi_k^{(k-2)}(0+) <0,$ we obtain $(-1)^{k-1} g^{(k-1)}(0+) =-\infty$ and contradict the induction hypothesis. Hence $(-1)^{k-1}\varphi_k^{(k-2)}(0+) \ge 0$ and by integration we get $(-1)^{k-1}\varphi_k\ge 0,$ which implies $(-1)^{k-1} g^{(k)}\ge 0$ and hence that $g'\in\cM_{k-1},$ too.

\endproof

\begin{remark}
\label{SkMk1}
{\em (a) As it turns out later, we will also have $(-1)^{k-1}\varphi_k^{(k-2)}(0+) = 0,$ because $b_k = 0$ with the notation of Equation (\ref{Psi1}) below. However, this fact seems more difficult to prove directly with the above induction argument.

\medskip

(b) We believe that $f\in\cS_k\Rightarrow g^{(i)}\in\cM_{k-i}$ for all $i =1, \ldots, k-1.$ This would give more precisions on the implication (\ref{Sti1}) $\Rightarrow$ (\ref{Sti0}).}
\end{remark}

\medskip

We can now prove the only if part of Theorem 1. Suppose that $f\in\cS_k$ and consider the following non-negative measure on $(0,\infty):$
$$\rho_k\; =\; (-1)^{k-1}\varphi_k^{(k-1)},$$
with the above notation. For every $x > 0,$ we have
\begin{equation}
\label{phi2}
\int_0^x \rho_k(dt)\; =\; (-1)^{k-1}\varphi_k^{(k-2)}(x)\, +\, (-1)^{k-2}\varphi_k^{(k-2)}(0+)\; < \; \infty
\end{equation}
by Proposition \ref{SkMk}. Suppose next that $\rho_k$ has a density, which is then the function
$$(-1)^{k-1}\varphi_k^{(k-1)}(x)\; =\; (-1)^{k-1} \Theta_k (g)(x)\; =\; (-1)^{k-1} \Theta_k ({\hat g})(x^{-1})\; =\; (-1)^{k-1}{\hat \varphi_k}^{(k-1)}(x^{-1}),$$
where the first two equalities come from Lemma \ref{Form}, and where we have used the symmetric notation ${\hat \varphi_n} = x^{2n-1} {\hat g}^{(n)}.$ Changing the variable implies
$$\int_x^\infty t^{-2} \rho_k(dt)\; =\; \int_0^{x^{-1}} \!\! {\hat \rho}_k(dt)$$
for every $x > 0,$ where ${\hat \rho}_k$ is the measure with density $(-1)^{k-1}{\hat \varphi_k}^{(k-1)}.$ By approximation, this equality remains true when $\rho_k$ is not necessarily absolutely continuous. Moreover, it is clear from the above argument applied to ${\hat g}$ that the right-hand side is also finite for every $x > 0.$ Thus we  have shown that $\rho_k$ integrates $1\wedge t^{-2},$ and hence the measure 
$$\nu_k(dt) \; =\; ((2k-2)! t)^{-1} \rho_k(dt),$$
integrates $t\wedge t^{-1}$ on $(0,\infty).$ Recalling the notation in (\ref{Msika}) for $\Psi_{k-1},$ we next observe that
\begin{eqnarray*}
((2k-2)! t)^{-1} \Psi_{k-1}^{(k)}(x,t) & = & ((2k-2)! t)^{-1} \frac{d^k}{dx^k} \lpa \psi_{k-1}(tx^{-1})\rpa\, \Un_{\{x\ge t\}}\\
& = &   \frac{(-1)^{k-1} }{(k-2)! x^{k+1}} \sum_{n=0}^{k-2} \lpa \!\!\!\begin{array}{c} k-2\\ n \end{array}\!\!\!\rpa (-tx^{-1})^{n}\Un_{\{x\ge t\}}.
\end{eqnarray*}
Therefore, we can compute
\begin{equation}
\label{Plug}
\int_0^\infty \Psi_{k-1}^{(k)}(x,t)\,\nu_k(dt)\; =\; \frac{(-1)^{k-1} }{(k-2)! x^{k+1}} \sum_{n=0}^{k-2} \lpa \!\!\!\begin{array}{c} k-2\\ n \end{array}\!\!\!\rpa (-x)^{-n}\lpa\int_0^x t^n\, \rho_k (dt)\rpa.
\end{equation}
Integration by parts with the help of Proposition \ref{SkMk} yields
$$\int_0^x t^n\, \rho_k (dt)\; =\; (-1)^{k-1}\,\sum_{i=0}^n \,(-1)^i \frac{n!\,x^{n-i}}{(n-i)!}\, \varphi_k^{(k-2-i)}(x)$$
for every $n\ge 1,$ the case $n=0$ being evaluated in the above (\ref{phi2}). Plugging all these expressions into (\ref{Plug}) and switching the two finite sums, we get
\begin{eqnarray}
\label{Psik}
\nonumber \int_0^\infty \Psi_{k-1}^{(k)}(x,t)\,\nu_k(dt) & = & \frac{1}{x^{2k-1}}\sum_{i=0}^{k-2}  \lpa \sum_{n=0}^i (-1)^n \binom{i}{n}\rpa \frac{x^i \varphi_k^{(i)}(x)}{i!}\; -\; \frac{\varphi_k^{(k-2)}(0+) }{(k-2)! x^{k+1}}\\
& =& \frac{\varphi_k(x)}{x^{2k-1}} \; -\; \frac{\varphi_k^{(k-2)}(0+) }{(k-2)! x^{k+1}}\; =\; g^{(k)} (x)\; -\; \frac{\varphi_k^{(k-2)}(0+) }{(k-2)! x^{k+1}}\cdot
\end{eqnarray}
Since $g'\in\cM_{k-1}$ by Proposition \ref{SkMk}, it is obvious that $g^{(i)}(\infty)$ exists and is finite for every $i =1,\ldots, k-1.$ Moreover we have
$$g^{(2)}(x^{-1})\; =\; -x^4 {\hat g}^{(2)}(x)\, +\, 2xg'(x^{-1})\; \to\; 0$$
as $x\to 0,$ by Proposition \ref{SkMk} applied to ${\hat g}.$ This yields $g^{(2)}(\infty) =0$ and clearly, we have $g^{(i)}(\infty) =0$ as well for every $i=3,\ldots,k-1.$ Moreover, since for every $i\ge 1$ and $x >1$ one has
$$\Psi_{k-1}^{(i)}(x,t)\; \le \; K_i\,x^{-1} (xt^{-1}\Un_{\{x\le t\}} + tx^{-1}\Un_{\{t<x\}})$$
for some finite constant $K_i,$ the integrability properties of $\nu_k$ show that
$$\int_0^\infty \Psi_{k-1}^{(i)}(x,t)\,\nu_k(dt)\; \rightarrow\; 0, \qquad x\to \infty$$
for every $i=1,\ldots,k-1.$ Hence, by monotone convergence, we can integrate $(k-1)$ times the identity (\ref{Psik})  from $x$ to $\infty$ and obtain
\begin{eqnarray}
\label{Psi1}
g' (x)\; =\; a_k\; +\; b_k\, x^{-2}\; +\; \int_0^\infty \Psi_{k-1}'(x,t)\,\nu_k(dt)
\end{eqnarray}
with $a_k = g'(\infty) \ge 0$ and
$$b_k\; =\; \frac{(-1)^{k-1}\varphi_k^{(k-2)}(0+) }{(k-2)!}\; \ge \; 0.$$
Integrating now (\ref{Psi1}) from $\varepsilon >0$ to $x,$ we get
$$g(x)\,-\, g(\varepsilon)\; =\; a_k(x-\varepsilon)\; +\; b_k(\varepsilon^{-1} -x^{-1})\; +\; \int_0^\infty (\Psi_{k-1}(x,t) -\Psi_{k-1}(\varepsilon,t))\,\nu_k(dt).$$
When $\varepsilon\downarrow 0,$ the left-hand side increases to $g(x)-g(0+) < \infty$ whereas (by Proposition \ref{IP1} and monotone convergence) the three terms on the right hand side increase to some limit which must be finite. This shows that $b_k = 0$ and, since $\Psi_{k-1}(0+,t) =0,$
$$g(x)\;=\; g(0+)\; +\; a_k \,x\; +\; \int_0^\infty \Psi_{k-1}(x,t)\,\nu_k(dt).$$
In particular, the measure $\nu_k(dt)$ must integrate $(1+ t)^{-1}.$ Dividing both sides by $x$ and setting
$$\mu_k(dt)\; =\; \nu_k(dt)\;+\; g(0+)\binom{2k}{k}^{-1}\delta_0(dt),$$
we have finally built a drift coefficient $a_k\ge 0$ and a non-negative measure $\mu_k(dt)$ on $[0,\infty)$ integrating $(1+t)^{-1},$ such that $f$ has the required representation
$$f(x)\; =\; a_k \; +\; \int_0^\infty \Phi_{k-1}(x,t)\,\mu_k(dt), \qquad x > 0.$$

\qed

\begin{remark}
\label{Inv}
{\em The above proof gives the following formula for the unique non-negative measure $\mu_k$ corresponding to $f\in\cS_k:$
$$\mu_k\; =\;  g(0+)\binom{2k}{k}^{-1}\delta_0\;+\; (-1)^{k-1}((2k-2)! t)^{-1}\varphi_k^{(k-1)},$$
with the notation $g=xf$ and $\varphi_k = x^{2k-1} g^{(k)}.$ This can be viewed as a Stieltjes inversion formula of finite type. This should be compared with Theorem 9 p. 345 and Theorem 10c p. 350 in \cite{W}, which give an analogous inversion formula for (\ref{SJ}), in the case when $\mu$ therein has a density.}
\end{remark}

\section{Proof of Theorem \ref{Bondlike}}

We begin with the following characterization of ${\widehat \cH_k}$ which has independent interest, and which will be used in both if and only if parts of the proof. For every non-negative differentiable function $f,$ let $\psi_f(x) = -x (\log f)'(x)$ and introduce for every fixed $u >0$ the function
$$\Delta_u(f)(w)\; =\; {\psi_f(uv)-\psi_f(uv^{-1}) \over v-v^{-1}}$$
for $v > 0,$ which is clearly a function of $w = v+v^{-1} \in[2,\infty)$ only.

\begin{lemma} \label{HMD} For $k\ge 2,$ one has
$$f\in {\widehat \cH_k}\; \Longleftrightarrow\; \Delta_u(f)\in\cM_{k-1}\;\;\,\forall\, u >0.$$
\end{lemma}

\proof Set $F_u(w) = \log(f(uv)f(uv^{-1}))$ for a given differentiable function $f.$ The crucial point is the following observation, which is obtained from the fact that $ dw/dv= v^{-1}(v-v^{-1})$ and the chain rule:
\begin{eqnarray}
\label{Dw}
\frac{d F_u}{dw}\; =\;  -\Delta_u (f).
\end{eqnarray}
This implies that if $\Delta_u (f)\in\cM_{k-1}$ for all $u >0,$ then $pF_u\in \cM_k$ for all $p,u >0.$ Moreover, it is easy to see from Fa\`a di Bruno's formula that $h\in \cM_k\Rightarrow e^h\in \cM_k$ for any given function $h.$ This concludes the if part of the lemma.

The only if part is analogous to the proof of Theorem 3.6 (iii) p. 19 in \cite{SSV}. Rewriting
$$f^p(uv)f^p(uv^{-1})\; =\; \sum_{n=0}^\infty \frac{p^n}{n!} \,(F_u(w))^n$$
and differentiating term by term, we see that if $f\in {\widehat \cH_k},$ then
$$(-1)^j \sum_{n=1}^\infty \frac{p^n}{n!} \,\lpa(F_u(w))^n\rpa^{(j)}\; \ge \; 0$$
for every $j=1,\ldots,k$ and every $p,u >0.$ Dividing by $p$ and letting $p\to 0,$ we get
$$(-1)^j F_u^{(j)}\; =\; (-1)^{j-1} \Delta_u (f)^{(j-1)}\; \ge \; 0$$
for every  $j=1,\ldots,k$ and $u > 0,$ as required.

\endproof

\begin{remark}
\label{DF}
{\em (a) The above proof shows that ${\widehat \cH_k} = \lacc f, \; \exists \,p_n\downarrow 0 \; \slash\; f^{p_n}\in\cH_k\racc.$ Below, we will see that there are examples of functions in $\cH_k$ such that $f^p\not\in \cH_k$ for some $p < 1.$ It is an open question whether $f^p\in \cH_k$ for every $p \ge 1$ as soon as $f\in\cH_k.$
\medskip

(b) It is easy to see from (\ref{Dw}) that if $f\in{\widehat \cH}_k$ for some $k\ge 2,$ then $\psi_f$ cannot take infinite values, so that necessarily $f(x)>0 $ for all $ x>0.$ This is in sharp contrast with the case $k=1.$ Observe in particular from Example 17.2.3 in \cite{LB97} that the $k-$monotone function
$$x\;\mapsto\; (t-x)_+^{k-1},$$
which is in $\cH_k$ for every $t > 0,$ cannot be in ${\widehat \cH}_k$ if $k\ge 2$ since it does not have full support.}
\end{remark}

We next state a crucial computational lemma.

\begin{lemma}
\label{DPsi}
For every $f$ regular enough and $k\ge 0, u >0,$ one has
$$ \lacc \Delta_u (f)^{(k)}(w)\racc_{w=2}\; = \; {k! \over (2k+1)!}\,
(u^{2k+1}\psi_f^{(k+1)}(u))^{(k)}.$$
\end{lemma}

\proof We will use the following polynomial identity
$$\sum_{k=0}^n \;z^k\; =\; \sum_{k=0}^{[n/2]}  \binom{n-k}{k}\, (-z)^k(1+z)^{n-2k},$$
where $[x]$ means the integer part of $x,$ which is an easy consequence of the Chu-Vandermonde identity. Setting now $x = v-1, y = v^{-1} -1$ and $X = w -2,$ we deduce
$$\frac{x^{n+1} -y^{n+1}}{x-y}\; =\; \sum_{k=0}^n \;x^ky^{n-k}\; =\; \sum_{k=0}^{[n/2]}  \binom{n-k}{k}\, (-xy)^k(x+y)^{n-2k}\; =\; \sum_{k=0}^{[n/2]}  \binom{n-k}{k}\, X^{n-k},$$
where we have used $x+y = -xy =X.$ Putting this together with a Taylor expansion of $y \mapsto \psi_f(uy)$ around $y =1,$ we obtain
\begin{eqnarray*}
\Delta_u (f)(w)& =& \sum_{n=0}^N \lpa \sum_{k=0}^{[n/2]}  \binom{n-k}{k}\, X^{n-k}\rpa
\frac{u^{n+1}\psi_f^{(n+1)}(u)}{(n+1)!}\; +\; O (X^{N+1})\\
& =& \sum_{n=0}^N \lpa \sum_{k=[(n+1)/2]}^n  \binom{k}{n-k}\, X^k\rpa
\frac{u^{n+1}\psi_f^{(n+1)}(u)}{(n+1)!}\; +\; O (X^{N+1})\\
& =& \sum_{k=0}^N \lpa \sum_{n=0}^k  \binom{k}{n}\, \frac{u^{n+k+1}\psi_f^{(n+k+1)}(u)}{(n+k+1)!}\rpa X^k
\; +\; O (X^{N+1})\\
& =& \sum_{k=0}^N \lpa {(u^{2k+1}\psi_f^{(k+1)}(u))^{(k)}\over (2k+1)!}\rpa X^k\; +\; O (X^{N+1}),\\
\end{eqnarray*}
where in the last equality we have used the Leibniz formula and the first equality in Lemma \ref{Form}. This  concludes the proof.

\endproof

\subsection{Proof of the only if part} We first consider the smooth case. Let $k\ge 2$ and suppose that $f\in\cC^{2k-1}\cap {\widehat \cH}_k.$  A combination of Lemmas \ref{HMD} and \ref{DPsi} shows that
$$(-1)^{n-1}(u^{2n-1}\psi_f^{(n)}(u))^{(n-1)}\; \ge \; 0$$
for every $u >0$ and $n=1,\ldots, k.$ A perusal of Section 2.2 shows that all this leads to the representation (\ref{Psi1}) for $g=\psi_f,$ that is
$$\psi_f'(x)\; =\;  a_k\; +\; b_k\, x^{-2}\; +\; \int_0^\infty \Psi_{k-1}'(x,t)\,\nu_k(dt)$$
for some $a_k, b_k \ge 0$ and $\nu_k$ a non-negative measure integrating $t\wedge t^{-1}.$ However, thinking e.g. of the function $f(x) = e^{-x^{-1}}$ which is $\HCM,$ we may have here $\psi_f(0+)=-\infty.$ We hence integrate this into
$$\psi_f(x)\; =\; c\; +\;  a_k\,x\; +\; \int_{[1,\infty)} \Psi_{k-1}(x,t)\,\nu_k(dt)\; -\; b_k\, x^{-1}\; +\; \int_{(0,1)} (\Psi_{k-1}(x,t)- P_{k-1}(0))\,\nu_k(dt)$$
for some constant $c \in\RR.$ It is not difficult to show that $\Psi_{k-1}(x,t)- P_{k-1}(0) = \Psi_{k-1}(x^{-1},t^{-1}),$ using $P_{k-1}(1) = -\tP_{k-1}(1).$  Hence, we can rewrite
\begin{eqnarray*}
\psi_f(x)& = & c\; +\;  \lpa a_k\,x\; +\; \int_{[1,\infty)} \Psi_{k-1}(x,t)\,\nu_k(dt)\rpa\; -\; \lpa b_k\, x^{-1}\; +\; \int_{(0,1)} \Psi_{k-1}(x^{-1},t^{-1})\,\nu_k(dt)\rpa\\
& = & c\; +\;  \lpa a_k\,x\; +\; \int_{[1,\infty)} \Psi_{k-1}(x,t)\,\nu_k(dt)\rpa\; -\; \lpa b_k\, x^{-1}\; +\; \int_{(1,\infty)} \Psi_{k-1}(x^{-1},t)\,{\hat \nu_k}(dt)\rpa
\end{eqnarray*}
where ${\hat \nu_k}$ is a non-negative measure on $(1,\infty)$ integrating $t^{-1}.$ Dividing by $x,$ we deduce
$$-(\log f)'(x)\; =\; c\,x^{-1}\; +\;  \lpa a_k\; +\; \int_{[1,\infty)} \Phi_{k-1}(x,t)\,\nu_k(dt)\rpa\; -\; x^{-2}\lpa b_k\; +\; \int_{(1,\infty)} \Phi_{k-1}(x^{-1},t)\,{\hat \nu_k}(dt)\rpa$$
which is, by Theorem 1 and setting $\beta = 1-c,$ the required representation of $f.$\\

Suppose last that $f\in{\widehat \cH}_k$ but is not necessarily $\cC^{2k-1}.$ Introduce the approximation
$$\psi_\varepsilon(x)\; =\; x^{-1} \int_0^\infty \psi_f(y)\,  \phi_\varepsilon (yx^{-1}) \,dy\; =\; \EE[\psi_f (xL_\varepsilon)]$$
where $h$ is a positive mollifier (for example $h(x) = \kappa e^{-(1-x^2)^{-1}}\Un_{\{\vert x\vert \le 1\}}$ where $\kappa$ is the normalizing constant), $h_\varepsilon(x) = \varepsilon^{-1} h(x\varepsilon^{-1})$ and
$$\phi_\varepsilon (x)\; =\; x^{-1} h_\varepsilon( \log x)$$
is the density of a random variable $L_\varepsilon$ with compact support $[e^{-\varepsilon},e^\varepsilon].$ The above integral is finite for every $x, \varepsilon >0$ since $f\in\pcH_2$ and is hence positive everywhere - see Remark \ref{DF} (b) above, so that $\psi_f$ is locally bounded on $(0,\infty).$ The same argument clearly shows that $\psi_\varepsilon$ is smooth. Setting
$$\Delta_u^\varepsilon (w)\; =\; {\psi_\varepsilon(uv)-\psi_\varepsilon(uv^{-1}) \over v-v^{-1}},$$
we get by Lemma \ref{HMD} and the linearity of the expectation that
$\Delta_u^\varepsilon\in\cM_{k-1}.$ Hence, the above argument shows that the representation
$$\psi_\varepsilon(x)\; = \; c_\varepsilon\; +\;  \lpa a_{k,\varepsilon}\,x\; +\; \int_1^\infty \Psi_{k-1}(x,t)\,\nu_{k, \varepsilon}(dt)\rpa\; -\; \lpa b_{k, \varepsilon}\, x^{-1}\; +\; \int_1^\infty \Psi_{k-1}(x^{-1},t)\,{\hat \nu_{k, \varepsilon}}(dt)\rpa$$
holds for every $\varepsilon > 0,$ for some non-negative measures $\nu_{k, \varepsilon}(dt)$ and ${\hat \nu_{k, \varepsilon}}(dt)$ integrating $t^{-1}$ at infinity. Since $\psi_\varepsilon\to \psi_f$ pointwise as $\varepsilon\downarrow 0,$ the conclusion follows from Helly's selection theorem.

\qed

\begin{remark}
\label{core}
{\em The above proof shows that $\cC^\infty\cap {\widehat \cH}_k$ is dense in ${\widehat \cH}_k$ for the pointwise topology. It is interesting to mention that $\cC^\infty\cap \cH_k$ is also dense in $\cH_k$ for the same topology. Indeed, if $f\in\cH_k,$ the approximation
$$f_\varepsilon(x)\; =\; \int_0^\infty f(y)\,  h_\varepsilon (xy^{-1}) \,y^{-1}dy,\qquad\quad
\mbox{where}\qquad\quad h_\varepsilon (x)\; =\; \frac{e^{-\frac{(\log x)^2}{2\varepsilon^2}}}{\sqrt{2\pi} \varepsilon x}$$
is the density of a log-normal distribution with variance parameter $\varepsilon^2$, is well-defined since $f\in\cH_1$ has the representation (\ref{H1}) and since $h_\varepsilon$ integrates any polynomial function at zero and infinity. Observe that by the change of variable $x\mapsto e^t,$ this amounts to the standard convolution approximation with a Gaussian kernel. In particular, one has $f_\varepsilon\in\cC^\infty$ and $f_\varepsilon\to f$ pointwise as $\varepsilon\downarrow 0.$ Finally, Example 17.2.5 and Property (iv) p. 302 in \cite{LB97} show that $f_\varepsilon\in\cH_k$ for all $\varepsilon > 0.$}
\end{remark}

\subsection{Proof of the if part}

Since $x^{\beta-1}\in\HCM$ and since the class ${\widehat \cH}_k$ is closed with respect to pointwise multiplication and to the transformation $f\mapsto{\tilde f}(x)=f(x^{-1}),$ it is enough to show that $f\in{\widehat \cH}_k$ whenever $-(\log f)'\in\cS_k.$ By  Lemma \ref{HMD}, Theorem 1, and monotone convergence, this amounts to showing that
\begin{eqnarray}
\label{CMK}
\Delta_{k,u}\;\in\;\cM_{k}
\end{eqnarray}
on $(2,\infty)$ for every $k\ge 1$ and $u >0$ where, recalling the definition of $\psi_k$ at the beginning of Section 2.1, we have set
$$\Delta_{k,u}(w)\;=\; {\psi_k(uv)-\psi_k(uv^{-1})\over v-v^{-1}}$$
and we have changed the parameter $k-1$ to $k$ for the simplicity of notation. \\

Suppose first $u =1.$ For every $v > 1,$ one has
\begin{eqnarray*}
\Delta_{k,1}(w) &=& \frac{1}{\sqrt{w^2-4}} \lpa \binom{2k}{k} \; +\; 2\,\sum_{n=1}^k \binom{2k}{k+n} (-v)^{-n}\rpa\\
& = & \frac{1}{\sqrt{w^2-4}} \lpa \binom{2k}{k} \;+ \;2\,\sum_{n=1}^k \binom{2k}{k+n} \lpa \frac{\sqrt{w^2-4} - w}{2}\rpa^n\rpa,
\end{eqnarray*}
and the same formula holds for $v < 1.$ Setting $x = (2-w)^{-1}<0,$ we next claim
that
$$\Delta_{k,1}(w)\; =\; x^{1-k}\lpa\sum_{n=0}^{k-1} \binom{2n}{n} \,x^n\; -\; (1-4x)^{-1/2}\rpa.$$
Indeed, both sides equal $1 - \sqrt{(w-2)(w+2)^{-1}}$ for $k=1$ and satisfy the recurrence relationship
$$u_{k+1} \, =\, x^{-1} u_k \, +\, \binom{2k}{k}, \quad k\ge 1.$$
If $w > 6$ viz. $4x \in (-1,0),$ we can transform the latter expression into
\begin{eqnarray}
\label{2F1}
\Delta_{k,1}(w) & = &  x^{1-k}\;\sum_{n\ge k} \binom{2n}{n} \,x^n \; =\; -x\,\binom{2k}{k}\,\pFq{2}{1}{1,k+1/2}{k+1}{4x}
\end{eqnarray}
where we have used the truncated binomial series formula - see e.g. Formula 2.8(9) p. 109 in \cite{EMOT} - for the second equality. Applying Euler's integral formula for the hypergeometric function - see e.g. Theorem 2.2.1 in \cite{AAR}, we finally get
$$\Delta_{k,1}(w)\; =\; \frac{(2k)!}{\sqrt{\pi} k! \Gamma (k+1/2)} \int_0^1 \lpa \frac{t^{k-1/2} (1-t)^{-1/2}}{w-2 +4t}\rpa \,dt,$$
a formula which remains true for $w > 2$ by analytic continuation. This shows that $\Delta_{k,1}\in\cCM$ on $(2,\infty)$ and readily implies (\ref{CMK}).\\

The proof in the case $u\neq 1,$ which is inspired by that of the main theorem in \cite{LB97}, is more subtle. Supposing first $w\in (2, u+u^{-1}),$ we have either $u > v > 1$ or $1>v > u^{-1}$ for $u >1,$ so that
$$\Delta_{k,u}(w) \; =\; {P_k(u^{-1}v^{-1})-P_k(u^{-1}v)\over v-v^{-1}}$$
for $u > 1.$ The same formula holds with $u$ replaced by $u^{-1}$ for $u < 1.$ This shows that $\Delta_{k,u}(w)$ is a polynomial of degree $k-1$ in $w.$ Moreover, it follows from Lemma \ref{DPsi}, the first equality in Lemma \ref{Form}, and (\ref{Sti2}) that for every $j = 0, \ldots, k-1,$
$$(-1)^j\Delta_{k,u}^{(j)}(2) \; \ge \; 0.$$
Hence, by Lemma \ref{Poly}, we will have $(-1)^j\Delta_{k,u}^{(j)}\ge 0$ on $(2, u+u^{-1})$ for every $j = 0, \ldots, k-1$ as soon as
\begin{eqnarray}
\label{upu}
(-1)^j\Delta_{k,u}^{(j)}(u+u^{-1}) \; \ge \; 0
\end{eqnarray}
for every $j = 0, \ldots, k-1.$ The latter is a consequence of the following simple  surprising formula, whose proof is postponed to the Appendix.

\begin{lemma}
\label{Rk}
For every $w\ge u+u^{-1}>2,$ one has
$$\Delta_{k,u}^{(k)}(w) \;  = \; {(-1)^{k}(2k)! (w-u-u^{-1})^k \over k!(w^2 - 4)^{k + 1/2}}\cdot
$$
\end{lemma}

We can now finish the proof of the if part of Theorem 2. It is plain by definition that $\Delta_{k,u}(w)$ and all its successive derivatives tend to zero as $w\to \infty.$ Hence, integrating successively from $w$ to $\infty$ the closed formula of Lemma \ref{Rk} shows
$$(-1)^j\Delta_{k,u}^{(j)}(w) \; \ge \; 0$$
for every $w\ge u +u^{-1}$ and $j = 0, \ldots, k.$ This implies (\ref{CMK}) on $[u+u^{-1}, \infty)$, and also on $(2,u+u^{-1})$ from the above considerations, since (\ref{upu}) holds true.

\qed

\begin{remark}
\label{Ini}
{\em The above proof shows the remarkable equivalence
$$\Delta_u(f)\in\cM_{k-1}\; \Longleftrightarrow\;(-1)^i\lacc\Delta_u(f)^{(i)}(w)\racc_{w=2}\,\ge\, 0\;\;\;\forall\, i=0, \ldots, k-1,$$
for every $u >0.$ In other words, the property $f\in\pcH_k$ is characterized only by the initial behaviour of the functions $\Delta_u(f), u >0.$}
\end{remark}

\section{Examples}

\label{Sec:Exa}

In this section we perform some explicit computations related to our main results, for some interesting classes of functions. We also define a family of positive self-decomposable distributions whose Laplace transforms are hyperbolically monotone of some finite order, but not in $\HCM.$

\subsection{Cauchy-type functions}
\label{Ex:Cauchy1}
We consider the functions
$$f_\alpha(x)\;= \;\frac{1}{1+2\cos(\pi\alpha) x+x^2}$$
with $\alpha\in[0,1).$ Such Cauchy-type functions appear in many situations, pure and applied. To give but one example, $f_\alpha(x)$ is a generating function of Tchebyshev polynomials of the second kind. More generally, $f_\alpha^p(x)$ is, for every $p >0,$ a generating function of Gegenbauer polynomials - see e.g. Formula (6.4.10) in \cite{AAR}.

\begin{proposition}
\label{Cauchy}
For every $k\ge 1,$ one has
$$f_\alpha\,\in\,{\widehat \cH_k}\;\Leftrightarrow\; 2\alpha k\le 1.$$
\end{proposition}

\proof It is clear that $f_0\in\HCM$ and we hence exclude the case $\alpha =0$ in the sequel. Computing
$$g_\alpha (x) \; =\; -(\log f_\alpha)'(x)\; =\; 2(\cos(\pi\alpha) + x)f_\alpha(x)\; =\;x^{-1}\lpa 2 - \frac{1}{1+xe^{\ii\pi\alpha}}- \frac{1}{1+xe^{-\ii\pi\alpha}} \rpa$$
implies, after some simplifications,
$$(-1)^{n-1}(x^n g_\a(x))^{(2n-1)}\; =\; (2n-1)!\,(f_\a(x))^{2n}\,\lpa e^{\ii n\pi\alpha}(x+e^{-\ii\pi\alpha})^{2n} + e^{-\ii n\pi\alpha}(x+e^{\ii\pi\alpha})^{2n} \rpa.$$
Hence, by the Corollary in Section 1, we have
$$f_\alpha\,\in\,{\widehat \cH_k} \;\Leftrightarrow\; e^{\ii n\pi\alpha}(x +e^{-\ii\pi\alpha})^{2n} + e^{-\ii n\pi\alpha}(x +e^{\ii\pi\alpha})^{2n} \;\ge\; 0\quad\forall\, n=1,\ldots,k\;\;\mbox{and}\; x>0,$$
which is equivalent to $\cos(n\pi\a) \ge 0$ for all $n=1,\ldots,k,$ and hence to $2\alpha k\le 1.$

\endproof

\begin{remark}
{\em (a) In this example, it is clear that $g_\a > 0$ with $(xg_\a)'(\infty) = 0$ and $(xg_\a)(0+) =0.$ By Remark \ref{Inv}, if $2\a k\le 1$ we have hence the representation
$$g_\a(x)\; =\; \int_0^\infty \Phi_{k-1}(x,t)\, \nu_{k,\a}(dt),$$
where $\nu_{k,\a}(dt)$ has density
\begin{eqnarray*}
f_{k,\a}(t)& =& (-1)^{k-1} ((2k-2)! t)^{-1}(t^{2k-1} (tg_\a(t))^{(k)})^{(k-1)} \\
& = & (-1)^{k-1} (2k-2)!^{-1}\, t^{k-1} (t^k g_\a(t))^{(2k-1)}\\
& = & (2k-1)\, t^{k-1} (f_\a(t))^{2k}\,( e^{\ii k\pi\alpha}(t+e^{-\ii\pi\alpha})^{2k} + e^{-\ii k\pi\alpha}(t+e^{\ii\pi\alpha})^{2k} ),\\
\end{eqnarray*}
which may rewritten as
$$f_{k,\a}(t)\; =\; \frac{(4k-2)t^{k-1}}{(1+2\cos(\pi\alpha) t+t^2)^{2k}}\times\lpa \sum_{n=0}^{2k} \binom{2k}{n} \cos((n-k)\pi\a)\, t^n\rpa.$$

\medskip

(b) It can be proved by elementary yet lengthy computations that
$$f_\a^p \in \cH_2 \; \Leftrightarrow\; \cos(\pi\a) \ge \frac{1}{\sqrt{2(p+1)}}\qquad\mbox{and}\qquad f_\a^p \in \cH_3 \; \Leftrightarrow\; \cos(\pi\a) \ge \sqrt{\frac{3}{2(p+2)}},$$
which shows that it might happen that $f_\a^p\in\cH_k$ and $f_\a^q\not\in\cH_k$ for some $q < p$ and $k=2,3.$ Observe also $f_\a^p \in \cH_1 \Leftrightarrow \cos(\pi\a) \ge 0.$ It is interesting to notice from Formula (6.4.11) in \cite{AAR} that 
$$\sqrt{\frac{3}{2(p+2)}}\qquad\mbox{resp.}\qquad \frac{1}{\sqrt{2(p+1)}}\qquad\mbox{resp.}\qquad 0$$ 
is the largest positive root of the Gegenbauer polynomial $C_3^\lambda(x)$ resp. $C_2^\lambda(x)$ resp. $C_1^\lambda(x)$ with $\lambda = p.$ Observe also that Proposition \ref{Cauchy} can be rephrased as
$$f_\alpha\,\in\,{\widehat \cH_k}\;\Leftrightarrow\; \cos(\pi\a) \ge \mu_{k,0+}$$
where $\mu_{k,0+} = \cos (\pi/2k)$ is the largest positive root of the Tchebyshev  polynomial of the first kind $T_k(x),$ which is itself a renormalized limit of $C_k^\lambda(x)$ as $\lambda \downarrow 0$ - see Formula (6.4.13) in \cite{AAR}. Recalling Remark \ref{DF} (a), it is hence natural to conjecture that
\begin{equation}
\label{Gegen}
f_\a^p \in \cH_k \; \Leftrightarrow\; \cos(\pi\a) \ge \mu_{k,p}
\end{equation}
where $\mu_{k,p}$ is the largest positive root of the Gegenbauer polynomial $C_k^\lambda(x)$ with $\lambda = p.$ This will be the matter of further research.

\medskip

(c) It is proved in \cite{TS14, TS15} that
$$f_\a (xy^{-1})\,\in\, {\rm TP}_k\; \Leftrightarrow\; \a\, <\, 1/k\;\;\;\mbox{or}\;\;\; \a\in\{1/k,\ldots, 1/2\}$$
for every $k\ge 2,$ which shows that one may have $f_\a (xy^{-1})\in {\rm TP}_k$ and $f_\a(x)\not\in{\widehat \cH_{k-1}}$ if $k\ge 3.$ This contrasts with the case $k=2$ where the two properties are equivalent, as recalled in the introduction. On the other hand, the above conjecture (\ref{Gegen}) reads
$$f_\a(x)\,\in\,\cH_{k-1}\; \Leftrightarrow\; \a\, \le\, 1/k$$
for $p=1,$ which shows that the two properties might have a similar characterization in this framework. Observe also that Conjecture 2.1 in \cite{TS14} relates the ${\rm TP}_\infty$ character of the kernel $f_\a^p (xy^{-1})$ to the above largest positive root $\mu_{k,p}.$

\medskip

(d) The 'imaginary extension' $f_{\ii\a}(x)$ admits the factorization
$$f_{\ii\alpha}(x)\;= \;\frac{1}{1+2\cosh(\pi\alpha) x+x^2}\; =\; \frac{1}{(e^{\pi\alpha} +x)(e^{-\pi\alpha} +x)}$$
and is hence clearly $\HCM$ for every real $\a.$ In particular, the kernel $f_{\ii\a} (xy^{-1})$ is always ${\rm TP}_2.$ On the other hand, the further total positivity properties of this interesting kernel do not seem easy to investigate at first sight.}
\end{remark}

\subsection{GIG-type densities}

We consider the probability density functions
$$f(x)\; =\; C x^{\beta -1}\exp(-a x^\alpha - b x^{-\a})$$
on $(0,\infty),$ with $a, b,\a > 0, \beta\in\RR$ and $C > 0$ is the normalizing constant. In the case $\a =1$, these densities correspond to the classical generalized inverse Gaussian (GIG) distributions and are also prototypes of functions in $\HCM,$ having the decomposition
$$f(x)\; =\; C x^{\beta -1}\, f_a (x) f_b(x^{-1})$$
with $-(\log f_a)' \equiv a$ and $-(\log f_b)' \equiv b$ both in $\cS.$ In the case $\a > 1,$ it is easy to see that these densities are in $\cH_1$ but not in $\cH_2.$ In the case $\a < 1$, it follows from general results that these densities are all in $\HCM$ - see \cite{LB92} p.60. More precisely, setting $a=b=1$ without loss of generality, we have the decomposition $f(x) = C x^{\beta -1} h_\a(x) h_\a(x^{-1})$ with
$$-(\log h_\a)'(x)\; =\; \a x^{\a-1}\; =\; \frac{\a \sin(\pi\a)}{\pi} \int_0^\infty \frac{t^{\a-1}}{t+x}\, dt\; \in\;\cS.$$
It is interesting to compare this formula with the finite type decompositions obtained in the present paper. For every $k\ge 2,$ Remark \ref{Inv} shows after some simplifications that we also have
$$\a x^{\a-1}\; =\; \int_0^\infty \Phi_{k-1}(x,t)\, f_{k,\a}(t)\,dt$$
with
$$f_{k,\a}(t)\; =\; (2k-2)!^{-1} \,\alpha^2\prod_{i=1}^{k-1}(i^2-\alpha^2)\; t^{\a-1}.$$
Observe that the Eulerian product formula for sines implies
$$\binom{2k-2}{k-1} f_{k,\a}(t)\; \rightarrow\; \lpa\frac{\a \sin(\pi\a)}{\pi}\rpa t^{\a-1}, \qquad k\to \infty,$$
as expected from (\ref{Stic}).

\subsection{Finite type Stieltjes functions}

In this paragraph we consider the functions $x\mapsto \Phi_k(x,1),$ which serve as building blocks for $\cS_{k+1},$ and hence belong to $\cS_n$ for every $n=2,\ldots, k.$ The proof of Proposition 1 implies after some computations that
\begin{equation}
\label{Phin}
\Phi_k(x,1)\; =\; \frac{1}{(2n)!}\,\int_0^\infty  \Phi_n(x,t)\, t^{-1}f_{n,k}(t\wedge t^{-1}) \,dt
\end{equation}
for every $n=1,\ldots, k-1,$ where $f_{n,k}(u) = (-1)^n u^{n+1} (u^n \psi_k(u))^{(2n+1)}.$ In particular, making $n=k-1$ shows that for $k\ge 1$ we have the recursive formula
$$ \Phi_k(x,1)\; =\; (2k-1)\,\int_0^\infty  \Phi_{k-1}(x,t)\,t^{-1}\,(t^k\wedge t^{-k}) \,dt.$$
Letting $k\to \infty$ in (\ref{Phin}) also implies, after some simplifications, the curious formula
$$\frac{1}{1+x}\;=\; (2n+1)\, \int_0^\infty  \Phi_n(x,t)\, \frac{t^n}{(1+t)^{2n+2}}\,dt$$
for every $n\ge 1,$ which is also a consequence of Remark \ref{Inv} and the first equality in Lemma \ref{Form}. Observe also that letting $n\to \infty,$ by (\ref{Stic}) and Stirling's formula this identity becomes a trivial one:
$$\frac{1}{1+x}\;=\; \int_0^\infty  \frac{1}{t+x} \;\delta_1(dt).$$

\subsection{A family of positive self-decomposable distributions}

We observe that the infinite positive measure $\nu_k (dx)$ with density
$$x^{-2}\Phi_k(x^{-1}, 1)\; =\; x^{-1} \psi_k(x^{-1})$$
is a L\'evy measure on $(0,\infty),$ that is it integrates $1\wedge x.$ Besides, since $x\mapsto \psi_k(x^{-1})$ is a non-increasing function by (\ref{Sti2}) with $n=1$, the infinitely divisible positive random variable $X_k$ with log-Laplace transform
$$-\log\EE[e^{-\lambda X_k}]\; =\; \varphi_k(\lambda)\; =\; \int_0^\infty (1- e^{-\lambda x})\, \nu_k(dx), \qquad \lambda >0,$$
is self-decomposable - see \cite{LB92} p. 18 or Proposition 5.15 in \cite{SSV}. On the other hand, the fact that $\psi_k$ is not smooth and hence not completely monotone prevents the Bernstein function $\varphi_k$ from being complete, so that a fortiori it is not Thorin-Bernstein either, with the terminology of Chapters 7 and 8 in \cite{SSV}. Hence, the function
$$\lambda\;\mapsto\; \EE[e^{-\lambda X_k}]$$
is not in $\HCM$ by Theorem 5.4.1 in \cite{LB92}. However, by Theorems 1 and 2 the latter function is in ${\widehat \cH_{k+1}},$ because
$$\varphi_k'(\lambda)\; =\; \int_0^\infty e^{-\lambda x}\, \psi_k(x^{-1})\, dx\; =\; \int_0^\infty \Phi_k(\lambda, t)\, e^{-t}\, dt.$$
Let us stress however that in general, non-increasing functions in $\pcH_k$ for some finite $k$ are not necessarily in $\cCM,$ as shows the example
$$\frac{1}{1 + 2\cos(\pi\a)\lambda + \lambda^2} \; =\; \frac{1}{\sin(\pi\a)} \int_0^\infty e^{-\lambda x}\, (e^{- \cos(\pi\a) x} \sin (\sin (\pi\a) x))\, dx,$$
which is for $\a\in (0, 1/2k)$ in $\pcH_k$  by Proposition \ref{Cauchy}, but not in $\cCM.$

\section{Appendix}
\label{Sec:App1}

In this section we prove Lemma \ref{Rk}, using an idea close to the argument of pp. 307-310 in \cite{LB97}. For $u>1,$ the condition $w \ge  u+u^{-1}$ implies either $v \ge u$ or $v\le u^{-1}.$ For $v\ge u$ we have
\begin{equation}
\label{Ignat} \Delta_{k,u}(w)\; =\; \frac{P_k(u^{-1}v^{-1})
+\tP_k(u v^{-1})}{v-v^{-1}}\; =\; {Q_k(v) \over v^{k} (v-v^{-1})}
\end{equation}
where, here and throughout, $Q_i(v)$ denotes any polynomial of
degree at most $i$ in $v.$ By the chain rule, using repeatedly $ dw/dv=
v^{-1}(v-v^{-1})$, we deduce
$$\Delta_{k,u}^{(k)}(w)\; =\;{Q_{3k}(v) \over v^{2k}(v-v^{-1})^{2k+1}}\cdot$$
On the other hand, since $P_k(x)$ is the polynomial part of the Laurent
series $(1-x)^k (1-x^{-1})^k,$ we have $P_k(x)\, +\, \tP_k(x^{-1})\, =\, (1-x)^k (1-x^{-1})^k \, =\, (-x)^{-k} (x-1)^{2k},$ which entails
$$\Delta_{k,u}(w)\; =\; (-1)^k\,\frac{u^k v^{-k}(v-u^{-1})^{2k}}{v-v^{-1}}\; +\; \frac{\tP_k(uv^{-1})
-\tP_k(uv)}{v-v^{-1}}\cdot$$
Observing that the second term on the right-hand side is a polynomial
of degree $k-1$ in $w,$ we deduce that the $k$-th derivative $\Delta_{k,u}^{(k)}(w)$ must factorize with $(v-u^{-1})^k$ and, by symmetry as a function of $w$, with $(v-u)^k$ as well. We hence obtain
\begin{equation}
\label{Sign} \Delta_{k,u}^{(k)}(w)\; =\;{(v-u)^k(v-u^{-1})^k
R_{k}(v) \over v^{2k}(v-v^{-1})^{2k+1}}
\end{equation}
where $R_{k}(v)$ is a polynomial of degree at most $k$ in $v.$ To get out
an expression for $R_k(v),$ we again use symmetry and consider the case $v\le u^{-1},$ where
\begin{equation}
\label{Times}
\Delta_{k,u}(w)\; =\; \frac{P_k(u^{-1}v) +\tP_k(u v)}{v^{-1}-v}\; =\; - {Q_k(v^{-1}) v^k\over (v-v^{-1})}
\end{equation}
where $Q_k$ is the same polynomial as in (\ref{Ignat}). It is clear that the sum of the two quantities in (\ref{Ignat}) and (\ref{Times}) is a polynomial of degree at most $k-1$ in $w,$ whose $k-$th derivative in $w$ is hence zero. This leads to
$$\Delta_{k,u}^{(k)}(w)\; =\;- {(v-u)^k(v-u^{-1})^k R_{k}(v) \over v^{2k}(v-v^{-1})^{2k+1}}$$
for $v\le u^{-1},$ with the same polynomial $R_k(v)$ as in (\ref{Sign}). Replacing now $v$ by $v^{-1}$ in (\ref{Sign}) and equating the two expressions, we
obtain $R_{k}(v) = v^{2k}  R_{k}(v^{-1}),$ which implies $R_k(v)
= C_{k,u} v^k$ for some non-zero constant $C_{k,u}$ depending on
$k,u$ only. Putting everything together, for $v \ge u$ we have
\begin{equation}
\label{Rad} \Delta_{k,u}^{(k)}(w)\;
=\;{C_{k,u}(v-u)^k(v-u^{-1})^k \over v^{k}(v-v^{-1})^{2k+1}}\;
=\;{C_{k,u} (w-u-u^{-1})^k \over (w^2-4)^{k+1/2}}
\end{equation}
and this formula is also clearly true for $v\le u^{-1}.$ In order to
identify the constant $C_{k,u},$ we let $w\rightarrow
\infty$ and obtain from (\ref{Ignat}) the behaviour
$$\Delta_{k,u}(w)\; \sim\; P_k(0)\, w^{-1}\; =\; \binom{2k}{k}\, w^{-1}.$$
Comparing with (\ref{Rad}), we finally obtain $C_{k,u} =(-1)^{k}
(2k)!/k!$ as required.

\qed

\begin{remark} {\em An alternative proof of Lemma 5 can be obtained with the help of the following exact formula:
$$\Delta_{k,u}(w) \; = \; -x\; \sum_{i=0}^k \binom{2k}{k+i} \,\pFq{2}{1}{i+1,k+1/2}{k+i+1}{4x}\, (xy)^i$$
for every $w\ge u+u^{-1}$ with the notations $x = (2-w)^{-1}$ and $y = u+u^{-1} -2.$ This identity, which extends (\ref{2F1}), turns out to be equivalent to some
bilinear formula of the Meixner type for the hypergeometric function - see \cite{EMOT} p. 84. Complete details have been written down and are available upon request.
Applying Euler's formula and making some simple hypergeometric transformations, we also get
$$\Delta_{k,u}^{(k)}(w) \; = \; {(-1)^{k}(2k)! (w-u-u^{-1})^k \over k!(w^2 - 4)^{k + 1/2}}\cdot$$
The unexpected point in Lemma \ref{Rk} is the simplicity of the expression of $\Delta_{k,u}^{(k)}(w),$ compared to that of $\Delta_{k,u}(w)$ and all its other derivatives in $w.$}
\end{remark}

\bigskip

\noindent
{\bf Acknowledgement.} The second author would like to thank Jean-Fran\c{c}ois Burnol for several discussions related to this paper.


\begin{thebibliography}{7}

\bibitem{AAR}
G.~E.~Andrews, R.~Askey and R.~Roy. {\em Special Functions.} Cambridge University Press, Cambridge, 1999.

\bibitem{BB}
A.~Behme and L.~Bondesson. A class of scale mixtures of gamma($k$)-distributions that are generalized gamma convolutions. To appear in {\em Bernoulli}.

\bibitem{Be}
C.~Berg. Quelques remarques sur le c\^one de Stieltjes. In: F.~Hirsch and G.~Mokobodzki (eds.) {\em S\'eminaire de Th\'eorie du Potentiel Paris, No. 5}. Lect. Notes Math. {\bf 814}, 70-79, 1980.

\bibitem{LB92}
L.~Bondesson. {\em Generalized Gamma Convolutions and Related Classes of Distributions and Densities.} Lect. Notes Stat. {\bf 76}, Springer-Verlag, New York, 1992.

\bibitem{LB97}
L.~Bondesson. On hyperbolically monotone densities. In: N.~L.~Johnson and N.~Balakrishnan (eds.) {\em Advances in the Theory and Practice of Statistics: A Volume in Honor of Samuel Kotz}, 299-313, Wiley, 1997.

\bibitem{EMOT}
A.~Erd\'elyi, W. Magnus, F.~Oberhettinger and F.~G.~Tricomi. {\em Higher Transcendental Functions Vol I.} McGraw-Hill, New-York, 1954.

\bibitem{H}
F.~Hirsch. Int\'egrales de r\'esolvantes et calcul symbolique. {\em Ann. Inst. Fourier} {\bf 22} (4), 239-264, 1972.

\bibitem{K}
S.~Karlin. {\em Total Positivity. Vol I.} Stanford University Press, Stanford, 1968.

\bibitem{P}
H.~L.~Pedersen. Pre Stieltjes functions. {\em Mediterr. J. Math.} {\bf 8} (1), 113-122, 2011.

\bibitem{SSV}
R.~L.~Schilling, R.~Song and Z.~Vondra\v{c}ek. {\em Bernstein Functions. Theory and Applications.} 2nd edition. De Gruyter, Berlin, 2012.

\bibitem{TS14}
T.~Simon. Total positivity of a Cauchy kernel. {\em J. Approx. Theory.} {\bf 184}, 238-258, 2014

\bibitem{TS15}
T.~Simon. Total positivity in stable semigroups. To appear in {\em Constructive Approximation.}

\bibitem{S}
A.~D.~Sokal. Real-variable characterization of generalized Stieltjes functions. {\em Expo. Math.} {\bf 28}, 179-185, 2010.

\bibitem{W0}
D.~V.~Widder. A classification of generating functions. {\em Trans. Amer. Math. Soc.} {\bf 39}, 244-298, 1936.

\bibitem{W1}
D.~V.~Widder. The Stieltjes transform. {\em Trans. Amer. Math. Soc.} {\bf 43}, 7-60, 1938.

\bibitem{W}
D.~V.~Widder. {\em The Laplace Transform.} Princeton University Press, Princeton, 1946.

\bibitem{Wi}
R.~E.~Williamson. Multiply monotone functions and their Laplace transforms. {\em Duke Math. J.} {\bf 23}, 189-207, 1956.

\end{thebibliography}
\end{document}